%change.tex:  In how many ways can I carry a total of n coins in my two pockets, and have the same amount in both pockets?`
%%a Plain TeX file by Shalosh B. Ekhad and Doron Zeilberger (x pages)

%begin macros

\baselineskip=14pt
\parskip=10pt

\font\eightrm=cmr8 
\font\eighttt=cmtt8
\magnification=\magstephalf

\def\1{{\overline{1}}}
\def\2{{\overline{2}}}
\parindent=0pt
\overfullrule=0in

\def\frac#1#2{{#1 \over #2}}
%\headline={\rm  \ifodd\pageno  \RightHead  \else  \LeftHead  \fi}
%\def\RightHead{\centerline{
%Title
%}}
%\def\LeftHead{ \centerline{Doron Zeilberger}}
%end macros
\centerline
{\bf In how many ways can I carry a total of n coins in my two pockets,} 
\centerline
{\bf and have the same amount in both pockets? }
\bigskip
\centerline
{\it Shalosh B. EKHAD and Doron ZEILBERGER}
\bigskip
\qquad\qquad\qquad\qquad {\it In fond memory of Gert Almkvist\footnote{$^*$}
{\eightrm Gert Almkvist was one of the most creative and original mathematicians that we have ever met.
He was known, among his friends, as ``the guy who generalized a mistake of Bourbaki"
[see \hfill\break 
{\eighttt http://sites.math.rutgers.edu/\~{}zeilberg/mamarim/mamarimhtml/gert.html}],
a master expositor (1989 Lester Ford award, joint with Bruce Berndt, and numerous articles in Swedish),
the co-inventor of the Almkvist-Zeilberger algorithm, and a great authority on Calabi-Yao differential
equations. In addition to his official affiliation with the University of Lund, he was
the founder of the Institute of Algebraic Meditation, and many of his papers used it as his affiliation.
}  (April 17, 1934- Nov. 24, 2018). 
}
\bigskip

{\bf Theorem 1}: Let $a(n)$ be the number of ways of having a total of $n$ coins in your two pockets
(each of them either a penny, a nickel, a dime, or a quarter), so that the amounts in
the pockets are identical, then
$$
\sum_{n=0}^{\infty} \, a(n) t^n \, = \, \frac{P(t)}{Q(t)} \quad , where
$$
$$
P(t) \, = \,
{t}^{54}+{t}^{53}+3\,{t}^{52}+4\,{t}^{51}+9\,{t}^{50}+15\,{t}^{49}+25\,{t}^{48}+37\,{t}^{47}+54\,{t}^{46}+76\,{t}^{45}+101\,{t}^{44}+128\,{t}^{43}
$$
$$
+158\,{t}^{42}+190\,{t}^{41}+226\,{t}^{40}+256\,{t}^{39}+290\,{t}^{38}+318\,{t}^{37}+353\,{t}^{36}+372\,{t}^{35}+394\,{t}^{34}+405\,{t}^{33}+425\,{t}^{32}
$$
$$
+431\,{t}^{31}+439\,{t}^{30}+438\,{t}^{29}+448\,{t}^{28}+448\,{t}^{27}+448\,{t}^{26}+438\,{t}^{25}+439\,{t}^{24}+431\,{t}^{23}+425\,{t}^{22}+405\,{t}^{21}
$$
$$
+394\,{t}^{20}+372\,{t}^{19}+353\,{t}^{18}+318\,{t}^{17}+290\,{t}^{16}+256\,{t}^{15}+226\,{t}^{14}+190\,{t}^{13}+158\,{t}^{12}+128\,{t}^{11}+101\,{t}^{10}
$$
$$
+76\,{t}^{9}+54\,{t}^{8}+37\,{t}^{7}+25\,{t}^{6}+15\,{t}^{5}+9\,{t}^{4}+4\,{t}^{3}+3\,{t}^{2}+t+1 \quad ,
$$
and
$$
Q(t) \, = \,
 \left( 1- t \right) ^{7} \left( 1+t \right) ^{5} \left( {t}^{2}+t+1 \right) ^{3} \left( {t}^{2}-t+1 \right) ^{2} \left( {t}^{12}+{t}^{11}+{t}^{10}+{t}^{9}+{
t}^{8}+{t}^{7}+{t}^{6}+{t}^{5}+{t}^{4}+{t}^{3}+{t}^{2}+t+1 \right)   \cdot
$$
$$
\left( {t}^{12}-{t}^{11}+{t}^{10}-{t}^{9}+{t}^{8}-{t}^{7}+{t}^{6}-{t}^{5}+{t}^{4}-{t}^{3
}+{t}^{2}-t+1 \right)  \left( {t}^{10}+{t}^{9}+{t}^{8}+{t}^{7}+{t}^{6}+{t}^{5}+{t}^{4}+{t}^{3}+{t}^{2}+t+1 \right)   \cdot
$$
$$
\left( {t}^{6}+{t}^{5}+{t}^{4}+{t}^{3}+{t}^{2}+t+1 \right) \quad .
$$

Furthermore $a(n)$ is a quasi-polynomial, and asymptotically,
$$
a(n) = \frac{5821}{311351040} \, n^6 \, + \, O(n^5) \quad .
$$
\vfill\eject
Finally
$$
a(10^{100})=
$$
1869593883482772371661260550149439038327927216816105704994593883482772371661260 \hfill\break
5501494390383279272213031310253532475754697976920199142421364643586865809088031 \hfill\break
3102535324757546979769201991424213646436436428797539908651019762130873241984353 \hfill\break
0954642065753176864287975399086510197621308732419843530954646363924141701919479 \hfill\break
6972574750352528130305908083685861463639241417019194796972574750352528130305908 \hfill\break
1017544132821910599688377466155243933021710799488577266355044132821910599688377 \hfill\break
4661552439330217107995247753069975292197514419736641958864181086403308625530847 \hfill\break
7530699752921975144197366419588641810864035 \quad .

The first $31$ terms are:
$$
1, 0, 4, 2, 12, 12, 34, 40, 85, 108, 190, 250, 394, 516, 762, 984, 1385, 1764,    2396, 
$$
$$
2998, 3966, 4886, 6316, 7684, 9739, 11706, 14594, 17358, 21320, 25134,   30470 \quad .
$$
As of Jan. 23, 2019, this sequence is not in the OEIS [Sl].

For analogous theorems where one can also have  a half-dollar coin, and a half-dollar coin as well as a dollar coin, see
the output files

{\tt http://sites.math.rutgers.edu/\~{}zeilberg/tokhniot/oEvenChange1b.txt}, and

{\tt http://sites.math.rutgers.edu/\~{}zeilberg/tokhniot/oEvenChange1c.txt} \quad .

{\bf How did we get this amazing theorem?}

let $A(n,m)$ be the number of ways of having $n$ coins in your two pockets (with denominations $1,5,10,25$) in such a way that
the difference between the amount in the left pocket and the amount in the right pocket is $m$ cents.
Then,  we have
$$
R(z,t):= \sum_{n=0}^{\infty} \sum_{m=-\infty}^{\infty} A(n,m)\,t^n z^m \quad,
$$
then
$$
R(z,t) \, = \, 
\frac{1}{(1-tz)(1-tz^5)(1-tz^{10})(1-tz^{25})} \cdot
\frac{1}{(1-t/z)(1-t/z^5)(1-t/z^{10})(1-t/z^{25})} \quad .
$$
This should be viewed as a {\it formal power series} in $t$ whose coefficients are {\it Laurent polynomials} in $z$, and we are
interested in extracting the coefficient of $z^0$. Now  you ask Maple to kindly convert the above rational function
into {\it partial fractions}, {\bf with respect to the variable ${\bf z}$}, getting something of the form
\vfill\eject
$$
\frac{P_1(z,t)}{1-tz} \, + \, \frac{P_2(z,t)}{1-tz^5} \, + \, \frac{P_3(z,t)}{1-tz^{10}} \, + \, \frac{P_4(z,t)}{1-tz^{25}}
\, + \,\frac{Q_1(z,t)}{z-t} \, + \, \frac{Q_2(z,t)}{z^5-t} \, + \, \frac{Q_3(z,t)}{z^{10}-t} \, + \, \frac{Q_4(z,t)}{z^{25}-t} \quad ,
$$
for some explicit  expressions in $z,t$, $P_1(z,t) \dots P_4(z,t)$,  $Q_1(z,t) \dots Q_4(z,t)$, that Maple finds for you.
These are {\bf rational functions} in $t$ but polynomials in $z$.

When we view them all as a formal power series in $t$, and take the coefficient of $z^0$, the $Q$'s do not contribute anything, so the
constant term, in $z$, is simply
$$
P_1(0,t) \, + \, P_2(0,t) \, + \, P_3(0,t) \, + \, P_4(0,t) \quad .
$$

This is implemented in the Maple package {\tt EvenChange.txt} by
procedure {\tt GfPAB(P,z,t,A,B)} that finds the coefficient of $z^0$ of
$$
\frac{P(z)}{\prod_{a \in A}(1-z^a \, t)\prod_{b \in B}(1-t/z^b)} \quad,
$$
for {\it any} polynomial $P$ of $z$ and {\it any} sets of positive integers $A$ and $B$
(so you can have different kinds of coins in each pocket, and also talk about the number of ways of doing it where the
difference between the amounts is not necessarily $0$).

The Maple package  {\tt EvenChange.txt} is available from the front of this article

{\tt http://sites.math.rutgers.edu/\~{}zeilberg/mamarim/mamarimhtml/change.html} \quad .

Since all the generating functions have  denominators whose roots are roots of unity, the sequence of interest itself, $a(n)$, is
a {\it quasi-polynomial}, albeit of a very large period. It is more efficient (still using partial fractions, this time with respect to $t$)
to express  it as a sum of quasi-polynomials of small periods. This is done via the procedure {\tt GFtoQPS} that is
lifted from the Maple package

{\tt http://sites.math.rutgers.edu/\~{}zeilberg/tokhniot/PARTITIONS} \quad ,

that accompanies [SiZ].

This is how we found so quickly $a(10^{100})$ and the leading asymptotics of $a(n)$ in Theorem 1.

{\bf Computing the generating functions  $\psi_n(t)$ dear to Gert Almkvist, Cayley, and Sylvester}

In his 1980 paper [A], Gert Almkvist  was interested in the sequence of rational functions $\{\psi_n(t)\}$ that he defined as
the constant term, in the variable $z$, of the rational function
$$
\frac{(1+z)^2}{2\,z \prod_{i=0}^{n} (1-t\,z^{n-2i})} \quad .
$$

Using procedure {\tt GfPAB} again, we got the following theorem. According to Almkvist [A], The cases $n=2,3,4$ are
due to Faa de Bruno and $n \leq 12$, except for $n=11$ are due to Sylvester and Franklin.

\vfill\eject

{\bf Theorem 2}:

$\bullet$ $n=2$:
$$
\psi_2(t)= {\frac {1}{ \left( 1 - t \right) ^{2} \left( t+1 \right) }} \quad,
$$
and its $n$-th coefficient, $a_2(n)$, is asymptotically
$$
a_2(n) \, = \, \frac{1}{2} n \, + \, O(1) \quad .
$$

The first $31$ terms (starting with $n=0$) are:
$$
1, 1, 2, 2, 3, 3, 4, 4, 5, 5, 6, 6, 7, 7, 8, 8, 9, 9, 10, 10, 11, 11, 12, 12,13, 13, 14, 14, 15, 15, 16 \quad .
$$
This is A4536 [{\tt http://oeis.org/A004526}] in  [Sl].

$\bullet$ $n=3$:
$$
\psi_3(t)= 
{\frac {{t}^{2}-t+1}{ \left( 1-t \right) ^{3} \left( t+1 \right)  \left( {t}^{2}+1 \right) }} \quad ,
$$
and its $n$-th coefficient, $a_3(n)$, is asymptotically
$$
a_3(n) \, = \, \frac{1}{8} n^2 \, + \, O(n) \quad .
$$

The first $31$ terms (starting with $n=0$) are:
$$
1, 1, 2, 3, 5, 6, 8, 10, 13, 15, 18, 21, 25, 28, 32, 36, 41, 45, 50, 55, 61,    66, 72, 78, 85, 91, 98, 105, 113, 120, 128 \quad .
$$
This is A1971 [{\tt http://oeis.org/A001971}] in  [Sl], that references [A].

$\bullet$ $n=4$:
$$
\psi_4(t)= 
{\frac {{t}^{2}-t+1}{ \left( 1- t \right) ^{3} \left( t+1 \right)  \left( {t}^{2}+1 \right) }} \quad ,
$$
and its $n$-th coefficient, $a_4(n)$, is asymptotically
$$
a_4(n) \, = \, \frac{1}{36} \, n^3 \, + \, O(n^2) \quad .
$$

The first $31$ terms (starting with $n=0$) are:
$$
1, 1, 3, 5, 8, 12, 18, 24, 33, 43, 55, 69, 86, 104, 126, 150, 177, 207, 241,    277, 318, 362, 410, 462, 519, 579, 645, 
$$
$$
715, 790, 870, 956 \quad .
$$

This is A1973 [{\tt http://oeis.org/A001973}] in  [Sl].

$\bullet$ $n=5$:
$$
\psi_5(t)= 
{\frac {{t}^{14}-{t}^{13}+2\,{t}^{12}+{t}^{11}+2\,{t}^{10}+3\,{t}^{9}+{t}^{8}+5\,{t}^{7}+{t}^{6}+3\,{t}^{5}+2\,{t}^{4}+{t}^{3}+2\,{t}^{2}-t+1}
{ \left( 1-t \right) ^{5} \left( t+1 \right) ^{3} \left( {t}^{2}+1 \right) ^{2} \left( {t}^{2}+t+1 \right)  \left( {t}^{2}-t+1 \right)  \left( {t}^{4}+1 \right) }}
\quad ,
$$
and its $n$-th coefficient, $a_5(n)$, is asymptotically
$$
a_5(n) \, = \, \frac{23}{4608} \, n^4 \, + \, O(n^3) 
\, = \, \frac{23}{2^9 \, 3^2} \, n^4 \, + \, O(n^3) 
\quad .
$$

The first $31$ terms (starting with $n=0$) are:
$$
1, 1, 3, 6, 12, 20, 32, 49, 73, 102, 141, 190, 252, 325, 414, 521, 649, 795, 967, 1165, 1394, 1651, 1944, 
$$
$$
2275, 2649, 3061, 3523, 4035, 4604, 5225, 5910 \quad .
$$
This is A1975 [{\tt http://oeis.org/A001975}] in  [Sl].

$\bullet$ $n=6$:
$$
\psi_6(t) = 
{\frac {{t}^{10}+{t}^{8}+3\,{t}^{7}+4\,{t}^{6}+4\,{t}^{5}+4\,{t}^{4}+3\,{t}^{3}+{t}^{2}+1}
{  \left( 1 -t \right)^{6} \left( {t}^{2}+1 \right)  \left( {t}^{4}+{t}^{3}+{t}^{2}+t+1
 \right)  \left( t+1 \right) ^{3} \left( {t}^{2}+t+1 \right) }} \quad , 
$$
and its $n$-th coefficient, $a_6(n)$, is asymptotically
$$
a_6(n) \, = \, \frac{11}{14400} \, n^5 \, + \, O(n^4) 
\, = \, \frac{11}{2^6\, 3^2 \, 5^2} \, n^5 \, + \, O(n^4) 
\quad .
$$

The first $31$ terms (starting with $n=0$) are:
$$
1, 1, 4, 8, 18, 32, 58, 94, 151, 227, 338, 480, 676, 920, 1242, 1636, 2137, 2739, 3486, 4370, 5444, 6698, 8196, 9926, 
$$
$$
11963, 14293, 17002, 20076, 23612, 27594, 32134 \quad .
$$
This is A1977 [{\tt http://oeis.org/A001977}] in  [Sl].

$\bullet$ $n=7$:
$$
\psi_7(t)= \frac{P_7(t)}{Q_7(t)} \quad,
$$
where
$$
P_7(t) =
{t}^{34}-{t}^{33}+3\,{t}^{32}+3\,{t}^{31}+7\,{t}^{30}+12\,{t}^{29}+16\,{t}^{28}+28\,{t}^{27}+33\,{t}^{26}+46\,{t}^{25}+56\,{t}^{24}+73\,{t}^{23}
$$
$$
+83\,{t}^{22}+90\,{t}^{21}+106\,{t}^{20}+109\,{t}^{19}+121\,{t}^{18}+110\,{t}^{17}+121\,{t}^{16}+109\,{t}^{15}+106\,{t}^{14}+90\,{t}^{13}+83\,{t}^{12}
$$
$$
+73\,{t}^{11}+56\,{t}^{10}+46\,{t}^{9}+33\,{t}^{8}+28\,{t}^{7}+16\,{t}^{6}+12\,{t}^{5}+7\,{t}^{4}+3\,{t}^{3}+3\,{t}^{2}-t+1 \quad ,
$$
and
$$
Q_7(t)=  \left( 1- t \right) ^{7} \left( t+1 \right) ^{5} \left( {t}^{2}+1 \right) ^{3} \left( {t}^{2}+t+1 \right) ^{2} \left( {t}^{2}-t+1 \right) ^{2}  \cdot
$$
$$
\left( {t}^
{4}+1 \right)  \left( {t}^{4}+{t}^{3}+{t}^{2}+t+1 \right)  \left( {t}^{4}-{t}^{3}+{t}^{2}-t+1 \right)  \left( {t}^{4}-{t}^{2}+1 \right) \quad ,
$$
and its $n$-th coefficient, $a_7(n)$, is, asymptotically
$$
a_7(n) \, =  \, \frac{841}{829440} \, n^6 \, + \, O(n^5) 
=  \, \frac{29^2}{2^{12} 3^4 5^2} \, n^6 \, + \, O(n^5)  \quad .
$$
The first $31$ terms (starting with $n=0$) are:
$$
1, 1, 4, 10, 24, 49, 94, 169, 289, 468, 734, 1117, 1656, 2385, 3370, 4672,    6375, 8550, 
$$
$$
11322, 14800, 19138, 24460, 30982, 38882, 48417, 59779, 73316,89291, 108108, 130053, 15564 \quad .
$$
This is A1979 [{\tt http://oeis.org/A001979}] in  [Sl].

For the cases $8 \leq n \leq 18$, see the output file

{\tt http://sites.math.rutgers.edu/\~{}zeilberg/tokhniot/oEvenChange2b.txt} \quad .

The case $n=8$ is A1981 [{\tt http://oeis.org/A001981}] in  [Sl].
As of Jan. 23, 2019, the cases $n=9$ and $n=10$ are not in the OEIS, and probably (we were too lazy to check) neither are the higher ones.

{\bf References}

[A] Gert Almkvist, {\it Invariants, mostly old ones}, Pacific J. Math. {\bf 86}(1980), 1-13. \hfill\break
{\tt https://projecteuclid.org/euclid.pjm/1102780612} \quad .

[SiZ] Andrew V. Sills and Doron Zeilberger,
{\it Formulae for the Number of Partitions of n into at most m parts(Using the Quasi-Polynomial Ansatz)},
Advances in Applied Mathematics {\bf 48} (2012), 640-645. \hfill\break
{\tt http://sites.math.rutgers.edu/\~{}zeilberg/mamarim/mamarimhtml/pmn.html} \quad .

[Sl] Neil A. J. Sloane, {\it The On-Line Encyclopedia of Integer Sequences}, {\tt  https://oeis.org/} \quad .

\bigskip
\hrule
\bigskip
Shalosh B. Ekhad, c/o D. Zeilberger, Department of Mathematics, Rutgers University (New Brunswick), Hill Center-Busch Campus, 110 Frelinghuysen
Rd., Piscataway, NJ 08854-8019, USA. \hfill\break
Email: {\tt ShaloshBEkhad at gmail dot com}   \quad .
\bigskip
Doron Zeilberger, Department of Mathematics, Rutgers University (New Brunswick), Hill Center-Busch Campus, 110 Frelinghuysen
Rd., Piscataway, NJ 08854-8019, USA. \hfill\break
Email: {\tt DoronZeil at gmail  dot com}   \quad .
\bigskip
\hrule
\bigskip
Exclusively published in the Personal Journal of Shalosh B.  Ekhad and Doron Zeilberger and arxiv.org \quad .
\bigskip
Written: Jan. 23, 2019.
\end